# Functional limit laws for the increments of the quantile process; with applications

## Vivian Viallon

*Laboratoire de Statistique Théorique et Appliquée, Université Paris VI, 175 rue du Chevaleret 75013 Paris FRANCE*
*e-mail:* `viallon@ccr.jussieu.fr`

**Abstract:** We establish a functional limit law of the logarithm for the increments of the normed quantile process based upon a random sample of size $n \to \infty$. We extend a limit law obtained by Deheuvels and Mason (12), showing that their results hold uniformly over the bandwidth $h$, restricted to vary in $[h'_n, h''_n]$, where $\{h'_n\}_{n \geq 1}$ and $\{h''_n\}_{n \geq 1}$ are appropriate non-random sequences. We treat the case where the sample observations follow possibly non-uniform distributions. As a consequence of our theorems, we provide uniform limit laws for nearest-neighbor density estimators, in the spirit of those given by Deheuvels and Mason (13) for kernel-type estimators.

## Contents

# 1. Introduction

## *1.1. Motivations*

The theory of empirical processes has been extensively investigated over the past decades. Several authors (see Deheuvels (7), Csörgő and Révész (4), Stute (27)) have underlined its relevance to the study of kernel nonparametric functional estimators, such as the Parzen-Rosenblatt density function estimator (refer to Parzen (22) and Rosenblatt (24)), or the Nadaraya-Watson regression function estimator (see, e.g., Nadaraya (21) and Watson (33)). Lately, these results have been significantly extended by, especially, Einmahl and Mason (14), (15), Mason (19), Deheuvels and Mason (13) and Varron (31). By combining empirical process arguments with combinatorial techniques initiated by Vapnik and Červonenkis ((30), (29)), these authors have established various versions of the uniform consistency of kernel type estimators, the latter being also shown to hold uniformly over the bandwidth parameter, with some restrictions. A direct statistical application of these results is given by the construction of limiting *certainty bands* for the density and the regression functions (see, e.g., Deheuvels and Mason (13)). A parallel field of study is that of quantile processes, for which a series of functional limit laws have been provided by Deheuvels and Mason (12) and Deheuvels (9). The motivation of the present study is that the recent refinements of the functional limit laws, initiated by Deheuvels, J.Einmahl, U.Einmahl and Mason, have been mostly written in the framework of the usual empirical (distribution) processes, and cover only in part the case of quantile processes. We will therefore orient our work towards bridging the remaining gaps in this theory. As a main result, we shall provide, in the sequel, a uniform in bandwidth functional limit law of the logarithm for the increments of the normed quantile processes. This result will be then applied to establish a uniform in bandwidth law of the logarithm for Nearest-Neighbor density function estimators. In this particular framework, we recall that such uniform in bandwidth results are particularly helpful to derive similar properties for estimators based on data-driven (and then random) bandwidths (see, for instance, (13)). We also mention that our results should be useful in the construction of nonparametric *goodness-of-fit* tests (see, e.g., Theorem 3.2 below). Likewise, some applications to the study of Lorenz process and score function estimators, in the spirit of Csörgő (6), may be derived from our results. We refer to Csörgő (6) for other examples of applications of quantile processes, which have some relevance to our work.

## *1.2. Notation and main Results*

Let $X_1, X_2, \cdots$ be independent and identically distributed [iid] random variables with distribution function $F(x) := \mathbb{P}(X_1 \leq x)$, for $x \in \mathbb{R}$, and quantile function $Q(t) := \inf\{x \geq 0 : F(x) \geq t\}$, for $0 \leq t \leq 1$. For $x \in \mathbb{R}$, we denote by $F_n(x) := n^{-1}\sharp\{X_i \leq x, 1 \leq i \leq n\}$ the empirical distribution function based upon $X_1, \cdots, X_n$ and by $Q_n(t) := \inf\{x \in \mathbb{R} : F_n(x) \geq t\}$, for $0 \leq t \leq 1$, the corresponding empirical quantile function. Here $\sharp E$ denotes the cardinality of $E$. We denote by $a_n(x) := n^{1/2}(F_n(x) - F(x))$, for $x \in \mathbb{R}$, the empirical process of order $n \geq 1$, and we let (see, e.g., Csörgő (6))

$$b_n(t) := n^{1/2}(Q_n(t) - Q(t))/q(t) \quad \text{for} \quad 0 < t < 1, \tag{1.1}$$

denote the normed quantile process. Here,

$$q(t) := \frac{d}{dt}Q(t) = \frac{1}{f(Q(t))}, \tag{1.2}$$

stands for the quantile density function (see, e.g., Parzen (23)), and

$$f(x) = \frac{d}{dx}F(x),$$

denotes the probability density function, both assumed to be properly defined and continuous in domains specified later on (see, e.g., the assumptions (F.1-2) below). In this paper, we are concerned with limit laws for the increments $\vartheta_n(\cdot)$ of the normed quantile process $b_n(\cdot)$, which are defined as follows. Given $0 < t < 1$ and a bandwidth $h \in (0, t \wedge (1-t))$, we set

$$\vartheta_n(t, h, s) := b_n(t + hs) - b_n(t), \text{ for } s \in [-1, 1]. \tag{1.3}$$

We will let $h > 0$ vary in such a way that $h'_n \leq h \leq h''_n$, where $\{h'_n\}_{n \geq 1}$ and $\{h''_n\}_{n \geq 1}$ are two sequences of positive constants such that $0 < h'_n \leq h''_n < 1$ and, for either choice of $h_n = h'_n$ or $h_n = h''_n$, the conditions (H.1-2-3) below are fulfilled by $\{h_n\}_{n \geq 1}$. Set, for $x \in \mathbb{R}^+$, $\log_2 x := \log_+ \log_+ x$, and $\log_+ x := \log(x \vee e)$. We assume that

(H.1)  $h_n \downarrow 0$ and $nh_n \uparrow \infty$ as $n \to \infty$;

(H.2)  $\log(1/h_n)/\log_2 n \to \infty$ as $n \to \infty$;

(H.3)  $nh_n/\log n \to \infty$ as $n \to \infty$.

In addition, we will say that the sequences $\{h'_n\}_{n \geq 1}$ and $\{h''_n\}_{n \geq 1}$, with $0 < h'_n \leq h''_n < 1$, for $n \geq 1$, fulfill the assumption (H.4) if at least one of the conditions (H.4)(i), (H.4)(ii) below are satisfied.

(H.4)(i) $\quad \dfrac{\sqrt{n}h'_n \log(1/h'_n)}{\log n \sqrt{\log_2 n}} \to \infty$ as $n \to \infty$;

(H.4)(ii) $\quad \dfrac{\sqrt{h''_n \log(1/h''_n)}}{h'_n \log(1/h'_n)} = o\left(\dfrac{\sqrt{n}}{\log n}\right)$ as $n \to \infty$.

**Remark 1.1.** (i) *Note that, under (H.1-2), the hypothesis (H.4)(i) is satisfied whenever (H.3′), introduced below, is satisfied.*

(H.3′) $\quad n^{1/2}h'_n/\log n \to \infty$ as $n \to \infty$.

*Obviously $(H.3')$ is a stronger condition than $(H.3)$.*

*(ii) For the particular choices $h'_n = n^{-r}$ and $h''_n = n^{-s}$, for $n \geq 1$, with $0 < s \leq r < 1$, the hypothesis $(H.4)(ii)$ is satisfied whenever $s \leq r < (1+s)/2$.*

*(iii) Our main results are established under the hypotheses $(H.1\text{-}2\text{-}3\text{-}4)$. However, it is unlikely that $(H.4)$ is necessary, and these results may still hold under the only hypotheses $(H.1\text{-}2\text{-}3)$ (see Remark 4.2 in the sequel).*

We now specify the range of $t$ in (1.3). We define the endpoints $u_1$ and $u_2$ of the random variable $X_1$ as follows

$$-\infty \leq u_1 := \inf\{x : F(x) > 0\} \quad \text{and} \quad \infty \geq u_2 := \sup\{x : F(x) < 1\}. \quad (1.4)$$

Further introduce $e_{n,h} := h + 25n^{-1}\log_2 n$, with $h \in [h'_n, h''_n]$. We will work on intervals of the form $[t_{1,n,h}, t_{2,n,h}]$ with $t_{1,n,h} = e_{n,h}$ and $t_{2,n,h} = 1 - e_{n,h}$.

We will assume that the following conditions hold.

$(F.1)$ F is twice continuously differentiable on $J := (u_1, u_2)$;

$(F.2)$ $F' = f$ is strictly positive on $J$;

$(F.3)$ $\sup_{u_1 < u < u_2} \left(F(u)(1 - F(u))|f'(u)|\right)/f^2(u) \leq \gamma$, for some $\gamma > 0$.

It is noteworthy that the inequality in $(F.3)$ is equivalent to

$$\sup_{0 < t < 1} t(1-t)\frac{|f'(Q(t))|}{f^2(Q(t))} \leq \gamma.$$

Some more notation is needed for the statement of our results. Denote by $\mathbb{S}_0$, the, so-called, Strassen set, a variant of which having been introduced by Strassen (25) in the framework of the law of the iterated logarithm for partial sums. Here, $\mathbb{S}_0$ is the unit ball of the reproducing kernel Hilbert space pertaining to the *two-sided standard Wiener process* $\{W(s) : |s| \leq 1\}$. The latter process is conveniently defined by setting

$$W(s) = \begin{cases} W_1(s) & \text{for } s \geq 0, \\ W_2(-s) & \text{for } s < 0, \end{cases}$$

where $W_1$ and $W_2$ are independent standard Wiener processes. We have, namely,

$$\mathbb{S}_0 = \left\{g \in AC(-1, 1), \quad g(0) = 0 \quad \text{and} \quad \int_{-1}^{1} \dot{g}(s)^2 ds \leq 1\right\},$$

where $AC(-1, 1)$ stands for the set of all absolutely continuous functions $g$ on $[-1, 1]$, with Lebesgue derivative $\dot{g}$. We denote by $\mathcal{B}(-1, 1)$ the set of all bounded functions on $[-1, 1]$ and set, for any $g \in \mathcal{B}(-1, 1)$, $\|g\| := \sup_{-1 \leq s \leq 1} |g(s)|$. Finally, we set, for any $\varepsilon > 0$, $\mathbb{S}_0^\varepsilon = \{h \in \mathcal{B}(-1, 1) : \inf_{g \in \mathbb{S}_0} \|g - h\| < \varepsilon\}$. Our main result may now be stated as follows.

**Theorem 1.1.** *Let $\{h'_n\}_{n\geq 1}$ and $\{h''_n\}_{n\geq 1}$ be two non-random sequences fulfilling the conditions (H.1-2-3-4), with $0 < h'_n \leq h''_n < 1$. Then, under (F.1-2-3), we have, almost surely,*

$$\lim_{n\to\infty} \sup_{h\in[h'_n,h''_n]} \left\{ \sup_{t_{1,n,h}\leq t\leq t_{2,n,h}} \left( \inf_{g\in\mathbb{S}_0} \left\| \frac{\vartheta_n(h,t;.)}{\sqrt{2h\log(1/h)}} - g \right\| \right) \right\} = 0. \qquad (1.5)$$

*Moreover, for each pair of constants $c_1, c_2$ with $0 \leq c_1 < c_2 \leq 1$, we have, almost surely,*

$$\forall g \in \mathbb{S}_0, \quad \lim_{n\to\infty} \sup_{h\in[h'_n,h''_n]} \left\{ \inf_{c_1\leq t\leq c_2} \left\| \frac{\vartheta_n(h,t;.)}{\sqrt{2h\log(1/h)}} - g \right\| \right\} = 0. \qquad (1.6)$$

**Remark 1.2.** *Extensions of Theorem 1.1 can be obtained by making sharper assumptions upon the distribution function $F$. In particular, $t_{1,n,h}$ and $t_{2,n,h}$ may be replaced in (1.5) by $h$ and $1-h$ respectively, when working under the additional assumptions (F.4-5) below.*

(F.4)  $A_1 = \lim_{u\downarrow u_1} f(u) < \infty, \quad A_2 = \lim_{u\uparrow u_2} f(u) < \infty;$

(F.5)  $\min(A_1, A_2) > 0.$

*It is noteworthy that assumptions (F.1-2-3-4-5) define particular tail monotone density functions with exponent $\gamma$, as defined by Parzen (23). Other less restrictive assumptions may lead to alternate extensions of Theorem 1.1. We refer the reader to Csörgő (6) for more details.*

The proof of Theorem 1.1 is postponed until Section 4. A rough outline of our arguments is as follows. First, we establish, in Section 2, a version of this theorem for iid uniform (0,1) random variables corresponding to the case where $f(x) = \mathbb{I}_{[0,1]}(x)$ and $q(t) = \mathbb{I}_{[0,1]}(t)$ (see, e.g., Theorem 2.1 below). Then, we make use of a continuity argument to treat the general framework. In Section 3 we present some statistical applications of our results. Section 4 is devoted to the proofs of our theorems. Finally, in the Appendix, we provide details on a technical fact used in Section 4.

## 2. The Uniform Case

Let $U_1, U_2, \cdots,$ be iid uniform $(0,1)$ random variables. In this context, we denote by $U_n(t) := n^{-1}\sharp\{U_i \leq t, 1 \leq i \leq n\}$, for $t \in \mathbb{R}$, the (right-continuous) empirical distribution function, and by $V_n(t) := \inf\{u \geq 0 : U_n(u) \geq t\}$, for $0 \leq t \leq 1$, with $V_n(t) = 0$ for $t \leq 0$ and $V_n(t) = V_n(1)$ for $t \geq 1$, the (left-continuous) empirical quantile function. Here, $\sharp E$ denotes the cardinality of $E$. We define the uniform *empirical process* based upon $U_1, \cdots, U_n$ by $\alpha_n(t) := n^{1/2}(U_n(t) - t)$ for $t \in \mathbb{R}$, and the corresponding uniform *quantile process*

$$\beta_n(t) := n^{1/2}(V_n(t) - t) \quad \text{for} \quad t \in \mathbb{R}. \qquad (2.1)$$

Note that, in this setup, the quantile process is equal to the normed quantile process, since, as was already pointed out, $q(t) = \mathbb{1}_{[0,1]}(t)$. Define further, for any $0 \leq t \leq 1$ and $h \in [0, t \wedge (1-t)]$, the increment functions

$$\xi_n(h,t;s) := \alpha_n(t+hs) - \alpha_n(t), \qquad (2.2)$$
$$\zeta_n(h,t;s) := \beta_n(t+hs) - \beta_n(t), \qquad (2.3)$$

for $-1 \leq s \leq 1$. The limiting behavior of the maximal oscillations of these processes has been extensively investigated in the literature (see, e.g., Stute (26), Deheuvels and Mason (11), (12), Deheuvels (9), Deheuvels and Einmahl (10), Mason (19), and the references therein).

We are now ready to state, in Theorem 2.1, a version of Theorem 1.1 for one-dimensional uniform $(0,1)$ random variables. It is noteworthy that Theorem 2.1 is somehow stronger than Theorem 1.1, in the sense that it holds on a larger interval (compare (1.5) and (2.4) below).

**Theorem 2.1.** *Let $\{h'_n\}_{n \geq 1}$ and $\{h''_n\}_{n \geq 1}$ be two non-random sequences fulfilling the conditions (H.1-2-3-4), with $0 < h'_n \leq h''_n < 1$. Then we have, almost surely,*

$$\lim_{n \to \infty} \sup_{h \in [h'_n, h''_n]} \left\{ \sup_{h \leq t \leq 1-h} \left( \inf_{g \in \mathbb{S}_0} \left\| \frac{\zeta_n(h,t;.)}{\sqrt{2h \log(1/h)}} - g \right\| \right) \right\} = 0. \qquad (2.4)$$

*Moreover, for any pair of constants $c_1, c_2$ with $0 \leq c_1 < c_2 \leq 1$, we have, almost surely,*

$$\forall g \in \mathbb{S}_0, \quad \lim_{n \to \infty} \sup_{h \in [h'_n, h''_n]} \left\{ \inf_{c_1 \leq t \leq c_2} \left\| \frac{\zeta_n(h,t;.)}{\sqrt{2h \log(1/h)}} - g \right\| \right\} = 0. \qquad (2.5)$$

The following corollary of Theorem 2.1 will be instrumental in the proof of Theorem 1.1, postponed until Section 4.2.

**Corollary 2.1.** *Under the assumptions of Theorem 2.1, we have, almost surely,*

$$\lim_{n \to \infty} \sup_{h \in [h'_n, h''_n]} \left\{ \sup_{h \leq t \leq 1-h} \left( \sup_{-1 \leq s \leq 1} \frac{|\zeta_n(h,t;s)|}{\sqrt{2h \log(1/h)}} \right) \right\} = 1. \qquad (2.6)$$

The proof of Theorem 2.1 is postponed until Section 4.1. An outline of our forthcoming arguments is as follows. First, we will establish, in Proposition 4.1 below, a version of this result for $\xi_n(h,t;\cdot)$. Proposition 4.1 will be shown to follow from Theorem 1.1 of Varron (31). Given this first result, Theorem 2.1 is straightforward under $(H.4)(i)$, and, to establish Theorem 2.1 under $(H.4)(ii)$, we will base the remainder of our proof on a uniform-in-bandwidth Bahadur-Kiefer-type representation of $\zeta_n$ in terms of $\xi_n$ (see, e.g., Bahadur (1), Kiefer (16), (17) and Deheuvels and Mason (12)). This representation is captured in Lemma 4.1 in the sequel.

## 3. Some Applications

### 3.1. The k-Spacings

In this sub-section, we provide some consequences of the just-given Theorems 1.1 and 2.1. The details of the corresponding proofs are postponed to Section 4.

We first consider the uniform case and denote by $U_{(1)} \leq \ldots \leq U_{(n)}$ the order statistics pertaining to a random sample $U_1, \ldots, U_n$ of iid uniform $(0,1)$ random variables. The *k-spacings* of $U_{(1)} \leq \ldots \leq U_{(n)}$ are then defined by

$$\Delta_{i,n}(k) := U_{(k+i)} - U_{(i)} \quad \text{for} \quad k = 1, \ldots, n \quad \text{and} \quad i = 0, \ldots, n+1-k, \quad (3.1)$$

where $U_{(0)} := 0$ and $U_{(n+1)} := 1$. For any integer $1 \leq d \leq n$, we set

$$\delta_n(d) := \max_{1 \leq k \leq d} \max_{0 \leq i \leq n+1-k} \left|\Delta_{i,n}(k) - k/n\right|. \quad (3.2)$$

The following Theorem 3.1 will be shown to follow from Theorem 2.1. Set $\lceil x \rceil \geq x > \lceil x \rceil - 1$ the ceiling function of $x \in \mathbb{R}$.

**Theorem 3.1.** *Let $\{h'_n\}_{n \geq 1}$ and $\{h''_n\}_{n \geq 1}$ be two non-random sequences fulfilling the conditions (H.1-2-3-4), with $0 < h'_n \leq h''_n < 1$. Then, we have, almost surely,*

$$\lim_{n \to \infty} \sup_{h \in [h'_n, h''_n]} \frac{\sqrt{n}\, \delta_n(\lceil nh \rceil)}{\sqrt{2h \log(1/h)}} = 1. \quad (3.3)$$

The proof of Theorem 3.1 is given in Section 4.4.

We now turn our attention to the case of possibly non-uniform random variables. Let $X_{(1)} \leq \ldots \leq X_{(n)}$ be the order statistics pertaining to a random sample $X_1, \ldots, X_n$, of iid random variables with common distribution function $F$, density function $f$ and quantile function $Q$. Let $e_n^{(1)} = 25n^{-1}\log_2 n$, and introduce the quantities $i_{1,n}$ and $i_{2,k,n}$ defined as follows

$$i_{1,n} := \min\{i : \frac{i}{n} \geq e_n^{(1)}\} \quad (3.4)$$

$$i_{2,k,n} := \max\{i : \frac{i+k}{n} \leq 1 - e_n^{(1)}\}. \quad (3.5)$$

Under the hypotheses (F.1-2-3), define the *k-spacings* of $X_{(1)} \leq \ldots \leq X_{(n)}$ by setting, for $k = 1, \ldots, n$ and $i = i_{1,n}, i_{1,n} + 1, \ldots, i_{2,k,n}$,

$$D_{i,n}(k) := X_{(k+i)} - X_{(i)}. \quad (3.6)$$

Note that the k-spacings $D_{i,n}(k)$ could obviously be defined for $k = 1, \ldots, n$ and $i = 0, \ldots, n - k$. However, the restriction we propose here is needed to derive Theorem 3.2, presented below, from Theorem 1.1.

For any integer $1 \leq d \leq n$, we set

$$d_n(d) := \max_{1 \leq k \leq d} \max_{i_{1,n} \leq i \leq i_{2,k,n}} f(X_{(i)}) \left| D_{i,n}(k) - \frac{k}{nf(X_{(i)})} \right|. \quad (3.7)$$

The following Theorem 3.2 provides a complement to Theorem 1 of Csörgő and Révész (5).

**Theorem 3.2.** *Let $\{h'_n\}_{n\geq 1}$ and $\{h''_n\}_{n\geq 1}$ be two non-random sequences fulfilling the conditions $(H.1\text{-}2\text{-}3\text{-}4)$, with $0 < h'_n \leq h''_n < 1$. Then, under $(F.1\text{-}2\text{-}3)$, we have, almost surely,*

$$\lim_{n\to\infty} \sup_{h\in[h'_n,h''_n]} \frac{\sqrt{n}\, d_n(\lceil nh\rceil)}{\sqrt{2h\log(1/h)}} = 1. \tag{3.8}$$

As mentioned above, Theorem 3.2 will be shown to be a consequence of Theorem 1.1, because of the well known relation $Q_n(t) = X_{(i)}$ for $(i-1)/n < t \leq i/n$. Its complete proof is postponed to Section 4.3.

### 3.2. A Law of the Logarithm for Nearest-Neighbor Density Estimators

In this sub-section, we show that the results of the Sub-section 3.1 imply a uniform-in-bandwidth law of the logarithm for a nearest-neighbor nonparametric density estimator. Let, as in Sub-section 3.1, $X_{(1)} \leq \ldots \leq X_{(n)}$ be the order statistics pertaining to a random sample $X_1, \ldots, X_n$ of iid variables with distribution function $F$, density function $F' = f$ and quantile function $Q$. Fix $0 < t_1 < t_2 < 1$ and define the random sequences

$$u_{1,n} := Q_n(t_1), \quad u_{2,n} := Q_n(t_2), \quad n \geq 1, \tag{3.9}$$

where $Q_n$ stands, as in Sub-section 1.2, for the empirical quantile function. Note that the fact that $q(t) = \frac{d}{dt}Q(t)$ exists and defines a positive and continuous function on $(0,1)$ is a consequence of (1.2) and (1.4), when combined with the assumptions $(F.1\text{-}2\text{-}3)$. Thus, for $n$ large enough, we have $u_1 < u_{1,n} < u_{2,n} < u_2$ almost surely, where $u_1$ and $u_2$ are defined in (1.4).

Further introduce $K$, an arbitrary kernel on $\mathbb{R}$, that is a measurable function integrating to one on $\mathbb{R}$, and denote by $\{k'_n\}_{n\geq 1}$ and $\{k''_n\}_{n\geq 1}$ two sequences such that $\{h'_n\}_{n\geq 1}$ and $\{h''_n\}_{n\geq 1}$ fulfill the conditions $(H.1\text{-}2\text{-}3\text{-}4)$, with $h'_n = k'_n/n$ and $h''_n = k''_n/n$, for $n \geq 1$. Select $k > 0$ such that $k \in [k'_n, k''_n]$. On the interval $[u_{1,n}, u_{2,n}]$, define the $k$ nearest-neighbor empirical density function, based upon the kernel $K$ and the sample $X_1, \ldots, X_n$, by

$$\widehat{f}_{n,k}(x) := \frac{1}{nR_k(x)} \sum_{i=1}^n K\Big(\frac{x-X_i}{R_k(x)}\Big), \tag{3.10}$$

where

$$R_k(x) := \inf\{r > 0, \text{ such that exactly } \lfloor k \rfloor \text{ elements of the sample}$$
$$X_1, \ldots, X_n \text{ are in } [x-r/2, x+r/2]\}.$$

This random function is often referred to as the *adaptative variable bandwidth of order k*. The following additional assumptions upon the kernel $K$ will be needed to state our result concerning nearest-neighbor density estimators.

(K.A)   $K$ is of bounded variation on $\mathbb{R}$.
(K.B)   $K$ is compactly supported.

**Theorem 3.3.** *Let $\{h'_n\}_{n \geq 1}$ and $\{h''_n\}_{n \geq 1}$ be two non-random sequences fulfilling the conditions (H.1-2-3-4), with $0 < h'_n \leq h''_n < 1$. Set, for $n \geq 1$, $k'_n = nh'_n$ and $k''_n = nh''_n$. Then, under (F.1-2-3) and (K.A-B), we have, almost surely,*

$$\lim_{n \to \infty} \sup_{x \in [u_{1,n}, u_{2,n}]} \sup_{k \in [k'_n, k''_n]} \frac{\sqrt{k}\, |\widehat{f}_{n,k}(x) - \mathbb{E}\widehat{f}_{n,k}(x)|}{\sqrt{2f^2(x) \log(n/k)}} = \left\{ \int_{\mathbb{R}} K^2(t) dt \right\}^{1/2}, \quad (3.11)$$

*where $u_{1,n}$ and $u_{2,n}$ are defined in (3.9).*

**Remark 3.1.** *Note that, in the setting of density estimation, the hypothesis $(H.4)(i)$ generally holds, since most of the bandwidth selection procedures lead to choices like $h'_n = c_1 n^{-\alpha}$, for a given $0 < c_1 < \infty$, with $\alpha \leq 1/2$.*

We will show, in the forthcoming Section 4.5, that Theorem 3.3 is a natural consequence of a combination of Corollary 2.1 of Varron (32), with Theorem 3.2. We mention that Theorem 1.1 of Varron (31) may be adapted likewise to obtain a similar result as that stated in Theorem 3.3, in the case of the usual Parzen-Rosenblatt (see, e.g., Parzen (22) and Rosenblatt (24)) kernel density estimator (see Facts 5.1 and 5.2 in the Appendix).

## 4. Proofs

### 4.1. Proof of Theorem 2.1

We first establish the analogue of Theorem 2.1 for $\xi_n$. Recall the definition (2.2) of $\xi_n$.

**Proposition 4.1.** *Let $\{h'_n\}_{n \geq 1}$ and $\{h''_n\}_{n \geq 1}$ be two non-random sequences fulfilling the conditions (H.1-2-3), with $0 < h'_n \leq h''_n < 1$. Then, with probability one,*

$$\lim_{n \to \infty} \sup_{h \in [h'_n, h''_n]} \left\{ \sup_{h \leq t \leq 1-h} \left( \inf_{g \in \mathbb{S}_0} \left\| \frac{\xi_n(h, t; .)}{\sqrt{2h \log(1/h)}} - g \right\| \right) \right\} = 0. \quad (4.1)$$

*Moreover, for any pair of constants $0 \leq c_1 < c_2 \leq 1$, we have, almost surely,*

$$\forall g \in \mathbb{S}_0, \quad \lim_{n \to \infty} \sup_{h \in [h'_n, h''_n]} \left\{ \inf_{c_1 \leq t \leq c_2} \left\| \frac{\xi_n(h, t; .)}{\sqrt{2h \log(1/h)}} - g \right\| \right\} = 0. \quad (4.2)$$

Note that the hypotheses (H.1-2-3) are sufficient in Proposition 4.1.

**Proof.** Set

$$\mathcal{G}_n := \{ \mathbb{I}_{[(t \wedge (t+v)), (t \vee (t+v))]}, \ h \in [h'_n, h''_n], \ h \leq t \leq 1-h, \ -h \leq v \leq h \}$$

and

$$\mathcal{G} := \{ \mathbb{I}_{[x,y]}, 0 \leq x \leq 1, 0 \leq y \leq 1 \}. \quad (4.3)$$

It is well known that the class of all closed intervals in $\mathbb{R}$ forms a $VC$ class (see, e.g., van der Vaart and Wellner (28)). Therefore, making use of the result of Exercise 9 on page 151 of van der Vaart and Wellner (28), it is readily shown that $\mathcal{G}$ constitutes a $VC$-subgraph class of functions (we refer to Section 2.6.2 in van der Vaart and Wellner (28) for the definitions of $VC$ classes of sets, and $VC$-subgraph classes of functions). Therefore $\mathcal{G}$ satisfies the entropy condition $(HK.III)(i)$ given in the Appendix, and since $\mathcal{G}_n \subset \mathcal{G}$, the equation (4.2) is a direct consequence of Theorem 1.1 of Varron (31) (which is recalled in Fact 5.1 in the Appendix for convenience). The same arguments readily show that, for each $0 < \lambda < 1/2$ and $\varepsilon > 0$, there exists almost surely an $n(\varepsilon)$ such that, for all $n \geq n(\varepsilon)$,

$$\{\xi_n(h,t;.) : h \in [h'_n, h''_n], t \in [\lambda, 1-\lambda]\} \subset \mathbb{S}_0^\varepsilon. \tag{4.4}$$

Our proof is completed by the observation that $\tilde{U}_i := [(U_i + \lambda) \text{ modulo } 1]$ and $U_i$ are identically distributed. Thus, for $n$ large enough, we have the distributional identities

$$\{\xi_n(h,t;.) : h \in [h'_n, h''_n], t \in [h, \lambda]\}_{n \geq 1}$$
$$\stackrel{d}{=} \{\xi_n(h,t;.) : h \in [h'_n, h''_n], t \in [\lambda, 2\lambda - h]\}_{n \geq 1},$$

and

$$\{\xi_n(h,t;.) : h \in [h'_n, h''_n], t \in [1-\lambda, 1-h]\}_{n \geq 1}$$
$$\stackrel{d}{=} \{\xi_n(h,t;.) : h \in [h'_n, h''_n], t \in [1-2\lambda, 1-\lambda-h]\}_{n \geq 1}.$$

By combining these statements with (4.4), we obtain, in turn, that for each $\varepsilon > 0$, there exists almost surely an $n(\varepsilon)$, such that for all $n \geq n(\varepsilon)$,

$$\{\xi_n(h,t;.) : h \in [h'_n, h''_n], t \in [h, \lambda]\} \subset \mathbb{S}_0^\varepsilon,$$

and

$$\{\xi_n(h,t;.) : h \in [h'_n, h''_n], t \in [1-\lambda, 1-h]\} \subset \mathbb{S}_0^\varepsilon.$$

The proof of proposition 4.1 in now complete.□

**Remark 4.1.** *We note that Corollary 3 of Mason (19), which is a d-variate extension of Theorem 3.1 of Deheuvels and Mason (12), can be established through these arguments.*

The completion of the proof of Theorem 2.1 under the hypothesis $(H.4)(i)$ is now straightforward. Indeed, Kiefer (16) showed that, almost surely as $n \to \infty$,

$$\sup_{0 \leq t \leq 1} |\alpha_n(t) + \beta_n(t)| = \mathcal{O}\Big(\frac{(\log n)^{1/2}(\log_2 n)^{1/4}}{n^{1/4}}\Big),$$

which implies that uniformly over $h \in [h'_n, h''_n]$, $h \leq t \leq 1-h$ and $|s| \leq 1$, we have almost surely as $n \to \infty$,

$$\frac{\zeta_n(h,t;s)}{\sqrt{2h\log(1/h)}} = \frac{-\xi_n(h,t;s)}{\sqrt{2h\log(1/h)}} + \mathcal{O}\Big(\frac{(\log n)^{1/2}(\log_2 n)^{1/4}}{n^{1/4}(h\log(1/h))^{1/2}}\Big),$$

and
$$\frac{\zeta_n(h,t;s)}{\sqrt{2h\log(1/h)}} = \frac{-\xi_n(h,t;s)}{\sqrt{2h\log(1/h)}} + \mathcal{O}\Big(\frac{(\log n)^{1/2}(\log_2 n)^{1/4}}{n^{1/4}(h\log(1/h))^{1/2}}\Big).$$

These two last results, when combined with Proposition 4.1, are enough to establish Theorem 2.1 under the hypothesis $(H.4)(i)$.

To prove Theorem 2.1 under the hypothesis $(H.4)(ii)$, we will work under the following notation. We will set $\|g\|_+ := \sup_{0 \le s \le 1} |g(s)|$ for the sup-norm of a function $g \in \mathcal{B}(0,1)$, in contrast with $\|g\| = \sup_{-1 \le s \le 1} |g(s)|$, used when $g \in \mathcal{B}(-1,1)$. To simplify matters, we will give proofs of our theorems with the formal replacement of $\|.\|$ by $\|.\|_+$. The technicalities needed to switch from $\|.\|_+$ to $\|.\|$ are straightforward, but lengthy, and will therefore be omitted.

The completion of the proof of Theorem 2.1 under the hypothesis $(H.4)(ii)$ will require a uniform-in-bandwidth Bahadur-Kiefer-type representation of $\zeta_n$ in terms of $\xi_n$. The following two lemmas are oriented towards the aim of establishing a representation of the kind. Our forthcoming results mimick that obtained by Deheuvels and Mason (12), in a slightly more general setup of varying bandwidths. We give the details of their proofs for the sake of completeness.

**Lemma 4.1.** *Let $\{h'_n\}_{n\ge 1}$ and $\{h''_n\}_{n\ge 1}$ be two sequences fulfilling the conditions $(H.1\text{-}2\text{-}3)$ and $(H.4)(ii)$. Then, for any $\lambda > 1$, we have, almost surely,*

$$\lim_{n\to\infty} \sup_{h\in[h'_n,h''_n]} \sup_{0\le t\le 1-\lambda h} \frac{\|\zeta_n(h,t;\cdot) + \xi_n(h, V_n(t);\cdot)\|_+}{\sqrt{2h\log(1/h)}} = 0, \qquad (4.5)$$

*where $V_n(\cdot)$ is still the empirical quantile function.*

**Proof.** Choose any $h \in [h'_n, h''_n]$, with $\lambda > 1$, $0 \le t \le 1 - \lambda h$ and $0 < s < 1$. Observe that

$$\zeta_n(h,t;s) + \xi_n(h, V_n(t); s) = \{\zeta_n(h,t;s) + \alpha_n(V_n(t+hs)) - \alpha_n(V_n(t))\} + \{\alpha_n(V_n(t) + hs) - \alpha_n(V_n(t+hs))\}. \qquad (4.6)$$

Making use of the easily proven fact that $|U_n(V_n(t)) - t| \le 1/n$ for any $n \ge 1$ and $0 \le t \le 1$, an application of the triangle inequality to the right-hand side of (4.6) establishes, in turn, that, uniformly over $h \in [h'_n, h''_n]$, $0 \le t \le 1 - \lambda h$ and $0 \le s \le 1$,

$$|\zeta_n(h,t;s) + \{\alpha_n(V_n(t+hs)) - \alpha_n(V_n(t))\}| \le 2n^{-1/2}. \qquad (4.7)$$

We invoke Theorem 1(III) in Mason (18) to obtain that, whenever $\{h''_n\}_{n\ge 1}$ satisfies $(H.1\text{-}2\text{-}3)$, we have

$$\lim_{n\to\infty} \sup_{0<h\le h''_n} \sup_{0\le t\le 1-h} \frac{\sqrt{n}|V_n(t+h) - V_n(t) - h|}{\sqrt{2h''_n\log(1/h''_n)}} = 1 \text{ a.s.}.$$

Thus, for any $\varepsilon > 0$, we have, almost surely for all $n$ sufficiently large, uniformly over $h \in [h'_n, h''_n]$, $0 \le t \le 1 - \lambda h$ and $0 \le s \le 1$,

$$|V_n(t+hs) - V_n(t) - hs| \le v_n := (1+\varepsilon)^2 \left(2n^{-1} h''_n \log(1/h''_n)\right)^{1/2}. \quad (4.8)$$

We next observe that, if the sequence $\{h''_n\}_{n \ge 1}$ satisfies, ultimately as $n \to \infty$, the assumptions $(H.1\text{-}2\text{-}3)$, then, such is also the case for $(v_n)_{n \ge 1}$ in (4.8). The first half of Proposition 4.1 implies therefore that

$$\limsup_{n \to \infty} \sup_{0 \le v \le v_n} \sup_{0 \le V_n(u) \le 1 - v_n} \frac{|\alpha_n(V_n(u)) - \alpha_n(V_n(u) + v)|}{\sqrt{2v_n \log(1/v_n)}} \le 1 \quad \text{a.s.} \quad (4.9)$$

(Note that $v$ in (4.9) plays the role of $hs$ in (4.1)).

Now, observe that, since $\{h''_n, n \ge 1\}$ satisfies $(H.1\text{-}2\text{-}3)$, we have

$$1/2 \le \liminf_{n \to \infty} \left(\log(1/v_n)\right)/\log n \le \limsup_{n \to \infty} \left(\log(1/v_n)\right)/\log n \le 1. \quad (4.10)$$

By combining (4.8) with (4.9) and (4.10), we get that, for all $n$ sufficiently large,

$$\sup_{h \in [h'_n, h''_n]} \sup_{0 \le t \le 1 - \lambda h} \sup_{0 \le s \le 1} n^{1/4} \left(h''_n \log(1/h''_n)\right)^{-1/4} (\log n)^{-1/2}$$
$$\times |\alpha_n(V_n(t+hs)) - \alpha_n(V_n(t) + hs)| \le 2^{3/4}(1+\varepsilon) \quad \text{a.s.}$$

Since $\varepsilon > 0$ can be chosen arbitrarily small in the above inequality, we see that, almost surely,

$$\limsup_{n \to \infty} \sup_{h \in [h'_n, h''_n]} \sup_{0 \le t \le 1 - \lambda h} \sup_{0 \le s \le 1} n^{1/4} \left(h''_n \log(1/h''_n)\right)^{-1/4} (\log n)^{-1/2}$$
$$\times |\alpha_n(V_n(t+hs)) - \alpha_n(V_n(t) + hs)| \le 2^{3/4}. \quad (4.11)$$

Under the hypothesis $(H.4)(ii)$, (4.11), when combined with (4.6) and (4.7), suffices to complete the proof of (4.5). □

**Remark 4.2.** *It is likely that the following result holds : for any sequences $\{h'_n\}_{n \ge 1}$ and $\{h''_n\}_{n \ge 1}$ fulfilling the hypotheses $(H.1\text{-}2\text{-}3)$,*

$$\lim_{n \to \infty} \sup_{h \in [h'_n, h''_n]} \sup_{0 \le t \le 1 - h} \sup_{0 \le s \le 1} \frac{|V_n(t+hs) - V_n(t) - hs|}{\sqrt{2h \log(1/h)}} \le 1 \; a.s. \; . \quad (4.12)$$

*If so, following the same lines as above, it could be shown that the result of Lemma 4.1, and, consequently, the results of Theorems 1.1, 3.1, 3.2 and 3.3, would still hold under the only hypotheses $(H.1\text{-}2\text{-}3)$.*

**Lemma 4.2.** *Let $\{h'_n\}_{n \ge 1}$ and $\{h''_n\}_{n \ge 1}$ be two sequences, assumed, each, to fulfill the conditions $(H.1\text{-}2\text{-}3)$ and $(H.4)(ii)$. Then, for any fixed $\lambda > 1$, we have almost surely*

$$\limsup_{n \to \infty} \sup_{h \in [h'_n, h''_n]} \sup_{1 - \lambda h \le t \le 1 - h} \frac{\|\zeta_n(h, t; \cdot) + \xi_n(h, 1-h; \cdot)\|_+}{\sqrt{2h \log(1/h)}} \le 2\sqrt{\lambda - 1}, \quad (4.13)$$

and for all $n$ sufficiently large, uniformly over $h \in [h'_n, h''_n]$,

$$V_n(1 - \lambda h) < 1 - h. \tag{4.14}$$

**Proof.** By setting, respectively, $t = h, s = -1$ in Proposition 4.1, and $t = 0, s = 1$ in Lemma 4.1, we obtain readily that, with probability one,

$$\limsup_{n \to \infty} \sup_{h \in [h'_n, h''_n]} \pm \beta_n(h)/\sqrt{2h \log(1/h)} \leq 1.$$

Thus, under ($H$.1-2-3), we have, almost surely,

$$\lim_{n \to \infty} \sup_{h \in [h'_n, h''_n]} \left(V_n(1 - \lambda h) - (1 - h)\right)/h = 1 - \lambda < 0. \tag{4.15}$$

Set $\omega_n(a) := \sup_{0 \leq t \leq 1-a} \{\|\xi_n(a, t; .)\|_+\}$ for $0 \leq a \leq 1$. From Proposition 4.1, we have $\sup_{h \in [h'_n, h''_n]}\{\omega_n(h)/\sqrt{2h \log(1/h)}\} \to 1$ almost surely as $n \to \infty$. By combining this result with (4.15), we obtain that

$$\limsup_{n \to \infty} \sup_{h \in [h'_n, h''_n]} \sup_{1 - \lambda h \leq t \leq 1-h} \frac{\|\xi_n(h, V_n(t); \cdot) - \xi_n(h, 1-h; \cdot)\|_+}{\sqrt{2h \log(1/h)}}$$
$$\leq \lim_{n \to \infty} \sup_{h \in [h'_n, h''_n]} \frac{2\omega_n((\lambda - 1)h)}{\sqrt{2h \log(1/h)}} = 2\sqrt{(\lambda - 1)} \quad \text{a.s.}$$

In view of (4.5), this last result readily yields (4.13), whereas (4.14) is a direct consequence of (4.15). $\square$

We are now ready to complete the proof of Theorem 2.1. To establish (2.4), we fix an $\varepsilon > 0$, and choose $\lambda = 1 + \frac{\varepsilon^2}{16}$. In view of (4.13), (4.14), Proposition 4.1 and (4.5), there exists almost surely an $N_\varepsilon < \infty$, such that, for all $n \geq N_\varepsilon$,

$$V_n(1 - \lambda h) < 1 - h \quad \text{uniformly over} \quad h \in [h'_n, h''_n],$$

and

$$\sup_{h \in [h'_n, h''_n]} \left\{ \sup_{0 \leq t \leq 1-h} \left( \inf_{g \in \mathbb{S}_0} \left\| \frac{\xi_n(h, t; \cdot)}{\sqrt{2h \log(1/h)}} - g \right\|_+ \right) \right\} < \varepsilon/2 \quad \text{a.s.},$$

$$\sup_{h \in [h'_n, h''_n]} \left\{ \sup_{0 \leq t \leq 1-\lambda h} \frac{\|\zeta_n(h, t; \cdot) + \xi_n(h, V_n(t); \cdot)\|_+}{\sqrt{2h \log(1/h)}} \right\} \leq \varepsilon/2 \quad \text{a.s.}, \tag{4.16}$$

$$\sup_{h \in [h'_n, h''_n]} \left\{ \sup_{1 - \lambda h \leq t \leq 1-h} \frac{\|\zeta_n(h, t; \cdot) + \xi_n(h, 1-h; \cdot)\|_+}{\sqrt{2h \log(1/h)}} \right\} \leq \varepsilon/2 \quad \text{a.s.}$$

Therefore, for all $n \geq N_\varepsilon$, we have, with probability one,

$$\sup_{h \in [h'_n, h''_n]} \left\{ \sup_{0 \leq t \leq 1-h} \left( \inf_{g \in \mathbb{S}_0} \left\| \frac{\zeta_n(h, t; \cdot)}{\sqrt{2h \log(1/h)}} - g \right\|_+ \right) \right\} < \varepsilon.$$

This last result suffices for the proof of the version of (2.4) obtained with $\|.\|_+$ replacing $\|.\|$. The proof of the version of (2.4) making use of the sup-norm $\|.\|$ follows along the same lines and is therefore omitted.

To establish (2.5), we first select an arbitrary $g \in \mathbb{S}_0$ and fix a $\varepsilon > 0$. By (4.2), there exist almost surely an $n_\varepsilon^{(4)}$ and a sequence $t_n^{(1)} \in (1/4, 3/4)$, $n = 1, 2, \ldots$, such that, for all $n \geq n_\varepsilon^{(4)}$,

$$\sup_{h \in [h_n', h_n'']} \| \xi_n(h, t_n^{(1)}; .)/\sqrt{2h \log(1/h)} - g \|_+ < \varepsilon/4. \tag{4.17}$$

Now set $t_n = U_n(t_n^{(1)})$ for $n \geq 1$. We have, uniformly over $h \in [h_n', h_n'']$,

$$\left\| \frac{\zeta_n(h, t_n^{(1)}; .)}{\sqrt{2h \log(1/h)}} + g \right\|_+ \leq \frac{\|\zeta_n(h, t_n; \cdot) + \xi_n(h, V_n(t_n); \cdot)\|_+}{\sqrt{2h \log(1/h)}}$$

$$+ \frac{\| -\xi_n(h, V_n(t_n); \cdot) + \xi_n(h, t_n^{(1)}; \cdot)\|_+}{\sqrt{2h \log(1/h)}}$$

$$+ \left\| -\frac{\xi_n(h, t_n^{(1)}; .)}{\sqrt{2h \log(1/h)}} + g \right\|_+. \tag{4.18}$$

The Glivenko-Cantelli theorem, when combined with the definition $t_n = U_n(t_n^{(1)})$ of $t_n$, with $t_n^{(1)} \in (1/4, 3/4)$, readily implies that, almost surely for all $n$ sufficiently large, $t_n \in (1/8, 7/8)$ and $V_n(t_n) \leq t_n^{(1)} < V_n(t_n + 1/n)$. This, in turn, entails that (see, e.g., Deheuvels (8))

$$\limsup_{n \to \infty} |V_n(t_n) - t_n^{(1)}|(n/\log n) = 1 \quad \text{a.s.} \tag{4.19}$$

Set $\rho_n := \log n/n$ for $n \geq 1$. By Theorem 1(I) of Mason, Shorack and Wellner (20), it follows that (see, e.g., (2.17) in Deheuvels and Mason (12))

$$\sup_{0 \leq t' \leq 1-\rho_n} \sup_{|t'-t''| \leq 2\rho_n} \frac{n}{\log n} \left[ U_n(t') - U_n(t'') \right] = \mathcal{O}(1) \text{ a.s.} .$$

By combining this last result with (4.19) and the fact that, under the hypotheses (H.1-2-3), $\log n/\sqrt{nh \log(1/h)} \to 0$ uniformly in $h \in [h_n', h_n'']$ as $n \to \infty$, we readily obtain that

$$\limsup_{n \to \infty} \sup_{h \in [h_n', h_n'']} \frac{\| -\xi_n(h, V_n(t_n); \cdot) + \xi_n(h, t_n^{(1)}; \cdot)\|_+}{\sqrt{2h \log(1/h)}} = 0 \quad \text{a.s.} . \tag{4.20}$$

By (4.16), (4.17) and (4.20), there exists almost surely an $N_\varepsilon' < \infty$, such that, for all $n \geq N_\varepsilon'$, we have $1/8 < t_n < 7/8$, $t_n < 1 - \lambda h_n''$ and

$$\sup_{h \in [h_n', h_n'']} \left\| \frac{\zeta_n(h, t_n^{(1)}; .)}{\sqrt{2h \log(1/h)}} - (-g) \right\|_+ < \varepsilon. \tag{4.21}$$

By choosing $\varepsilon > 0$ arbitrarily small in (4.21), we obtain readily (2.5) with $c_1 = 1/4$ and $c_2 = 3/4$. The extension to arbitrary $0 < c_1 \leq c_2 < 1$ follows from the same lines and is then omitted for the sake of conciseness. This last result completes the proof of the version of Theorem 2.1 pertaining to the sup-norm $\|.\|_+$. The extension to the case of $\|.\|$ is straightforward, and hence also omitted. □

### 4.2. Proof of Theorem 1.1

Recall that $Q_n(t) = Q(V_n(t))$ for $0 < t < 1$. Keep in mind that the fact that $q(t) = \frac{d}{dt}Q(t)$ exists, and defines a positive and continuous function on $(0,1)$, is a consequence of (1.2) and (1.4), when combined with the assumptions ($F$.1-2-3). In view of the definitions (1.1)-(2.1) of $b_n$ and $\beta_n$, this, when combined with Taylor's formula, entails, almost surely for all $n$ sufficiently large, the existence of a $\theta_{t,n} \in (t \wedge V_n(t), t \vee V_n(t))$ for each $t \in (0,1)$, such that

$$b_n(t) = \frac{n^{1/2}(Q(V_n(t)) - Q(t))}{q(t)} = \beta_n(t)\frac{q(\theta_{t,n})}{q(t)}. \tag{4.22}$$

Thanks to this relation, Csörgő and Révész (3) showed that, under the assumptions ($F$.1-2-3) (see also Theorem 3.2.1 in (6))

$$\sup_{t \in [e_n^{(1)}, 1-e_n^{(1)}]} |b_n(t) - \beta_n(t)| = \mathcal{O}\Big(\frac{\log_2 n}{n^{1/2}}\Big) \text{ almost surely,} \tag{4.23}$$

with $e_n^{(1)} = 25n^{-1}\log_2 n$. The result of Theorem 1.1 directly follows by combining Theorem 2.1 with (4.23).

### 4.3. Proof of Theorem 3.2

Recall the definition (1.1) of $b_n(t) = n^{1/2}(Q_n(t) - Q(t))/q(t)$, where $Q_n(t) = X_{(i)}$ for $(i-1)/n < t \leq i/n$. Thus, for any $0 < s, t < 1$, setting $i = \lceil nt \rceil$ and $j = \lceil ns \rceil$, where $\lceil x \rceil \geq x > \lceil x \rceil - 1$ is the ceiling function, we have $Q_n(s) = X_{(j)}$, $Q_n(t) = X_{(i)}$, and

$$b_n(s) - b_n(t) = n^{1/2}\Big\{\frac{X_{(j)} - Q(s)}{q(s)} - \frac{X_{(i)} - Q(t)}{q(t)}\Big\}$$
$$= \frac{n^{1/2}}{q(t)}\Big\{X_{(j)} - X_{(i)} - (Q(s) - Q(t)) + (X_{(j)} - Q(s))\Big(\frac{q(t)}{q(s)} - 1\Big)\Big\}$$
$$= \frac{n^{1/2}}{q(t)}\Big\{\Big(X_{(j)} - X_{(i)} - (Q(\frac{j}{n}) - Q(\frac{i}{n}))\Big) + \Big((Q(\frac{j}{n}) - Q(\frac{i}{n}))$$
$$- (Q(s) - Q(t))\Big) + (Q_n(s) - Q(s))\Big(\frac{q(t)}{q(s)} - 1\Big)\Big\}.$$

This allows us to write

$$b_n(s) - b_n(t) = \frac{n^{1/2}}{q(t)}\left[\left(X_{(j)} - X_{(i)} - Q\left(\frac{j}{n}\right) - Q\left(\frac{i}{n}\right)\right) + \varepsilon_n(t,s)\right], \quad (4.24)$$

where

$$\begin{aligned}
\varepsilon_n(t,s) &:= \varepsilon_{1,n}(t,s) - \varepsilon_{2,n}(t,s), & (4.25) \\
\varepsilon_{1,n}(t,s) &:= (Q_n(s) - Q(s))\left(\frac{q(t)}{q(s)} - 1\right), \\
\varepsilon_{2,n}(t,s) &:= Q(s) - Q(t) - \left(Q\left(\frac{j}{n}\right) - Q\left(\frac{i}{n}\right)\right).
\end{aligned}$$

Since $(i-1)/n < t \leq i/n$ and $(j-1)/n < s \leq j/n$, the assumptions (F.1–2) readily imply that,

$$|Q(s) - Q\left(\frac{j}{n}\right)| \leq |s - \frac{j}{n}|q(s_1) \leq \frac{q(s_1)}{n} \quad \text{for some } s_1 \in ((j-1)/n, j/n),$$

and, similarly,

$$|Q(t) - Q\left(\frac{i}{n}\right)| \leq |t - \frac{i}{n}|q(t_1) \leq \frac{q(t_1)}{n} \quad \text{for some } t_1 \in ((i-1)/n, i/n).$$

It follows therefore that we have, uniformly over $h \in [h'_n, h''_n]$, ultimately as $n \to \infty$,

$$\sup_{e_n^{(1)} \leq t \leq t_{2,n,h}} \sup_{0 \leq u \leq 1} |\varepsilon_{2,n}(t, t+hu)| = \mathcal{O}\left(\frac{1}{n}\right) \text{ a.s. .} \quad (4.26)$$

To evaluate $\varepsilon_{1,n}(t, t+hu)$, we will make use of the following fact, due to Csörgő and Révész (3) (see also Lemma 1.4.1 in Csörgő (6)).

**Fact 4.1.** *Under (F1-2-3-4), we have, for every $y_1, y_2 \in (0,1)$,*

$$\frac{q(y_2)}{q(y_1)} = \frac{f(Q(y_1))}{f(Q(y_2))} \leq \left\{\frac{(y_1 \vee y_2)}{(y_1 \wedge y_2)} \times \frac{1 - (y_1 \wedge y_2)}{1 - (y_1 \vee y_2)}\right\}^\gamma,$$

*where $\gamma > 0$ is as in (F.3).*

Under (F.3), it follows readily from Fact 4.1 that, for $u > 0$ and as $h \to 0$,

$$\frac{q(t)}{q(t+hu)} - 1 = \mathcal{O}(h^\gamma). \quad (4.27)$$

Similarly, for $u < 0$, one can show that $\frac{q(t)}{q(t+hu)} - 1 = \mathcal{O}(h^\gamma)$ as $h \to 0$. In addition, by the Chung (2) law of the iterated logarithm [LIL] applied to the sup-norm of $\beta_n$, we observe that, almost surely as $n \to \infty$,

$$n^{-1/2} \sup_{0 \leq t \leq 1} |\beta_n(t)| = \mathcal{O}\left(\sqrt{\frac{\log_2 n}{n}}\right). \quad (4.28)$$

We now combine (4.27) with (4.22) and (4.28). We so obtain that we have, with probability 1, uniformly over $h \in [h'_n, h''_n]$ and ultimately as $n \to \infty$,

$$\sup_{e_n^{(1)} \leq t \leq t_{2,n,h}} \sup_{0 \leq u \leq 1} |\varepsilon_{1,n}(t, t+hu)| = \mathcal{O}\Big(h^\gamma \sqrt{\frac{\log_2 n}{n}}\Big). \tag{4.29}$$

Recall that $(i-1)/n < t \leq i/n$ in (4.24). A Taylor expansion based upon (F.1), together with the Chung (2) LIL, as stated in (4.28), and (4.22), shows readily that, almost surely as $n \to \infty$,

$$\begin{aligned}
|f(Q(t)) - f(X_{(i)})| &= |f(Q(t)) - f(Q_n(t))| = \mathcal{O}\Big(|Q(t) - Q_n(t)|\Big) \\
&= \mathcal{O}\Big(|Q(t) - Q(\tfrac{i}{n})|\Big) + \mathcal{O}\Big(|Q(\tfrac{i}{n}) - Q_n(\tfrac{i}{n})|\Big) \\
&= \mathcal{O}\Big(\sqrt{\tfrac{\log_2 n}{n}}\Big) + \mathcal{O}\Big(\tfrac{1}{n}\Big) = \mathcal{O}\Big(\sqrt{\tfrac{\log_2 n}{n}}\Big). \quad (4.30)
\end{aligned}$$

Set

$$d_n^{(1)}(d) := \max_{1 \leq k \leq d} \max_{i_{1,n} \leq i \leq i_{2,k,n}} f(X_{(i)}) \big| D_{i,n}(k) - \big(Q(\tfrac{i+k}{n}) - Q(\tfrac{i}{n})\big)\big|,$$

where we recall that

$$\begin{aligned}
i_{1,n} &= \min\{i : \tfrac{i}{n} \geq e_n^{(1)}\} \\
i_{2,n} &= \max\{i : \tfrac{i+k}{n} \leq 1 - e_n^{(1)}\},
\end{aligned}$$

with $e_n^{(1)} = 25 n^{-1} \log_2 n$. Observe that, for all $h \in [h'_n, h''_n]$ and every $i_{1,\lceil nh \rceil, n} \leq i \leq i_{2,n}$,

$$\forall t \in \Big(\frac{i-1}{n}, \frac{i}{n}\Big], \quad Q_n(t) = X_{(i)}, \tag{4.31}$$

$$\exists t \in \Big(\frac{i-1}{n}, \frac{i}{n}\Big], \quad Q_n(t+h) = X_{(i+\lceil nh \rceil)}, \tag{4.32}$$

$$\not\exists t \in \Big(\frac{i-1}{n}, \frac{i}{n}\Big], \quad Q_n(t+h) = X_{(i+\lceil nh \rceil+1)}. \tag{4.33}$$

Now, looking carefully at the arguments used in the proof of Theorems 1.1 and 2.1, it is straightforward that the following results hold. Under the assumptions of Theorem 1.1, we have, almost surely,

$$\lim_{n \to \infty} \sup_{h \in [h'_n, h''_n]} \Big\{ \sup_{e_n^{(1)} \leq t \leq t_{2,n,h}} \Big( \sup_{0 \leq s \leq 1} \frac{|\vartheta_n(h, t; s)|}{\sqrt{2h \log(1/h)}} \Big) \Big\} = 1. \tag{4.34}$$

In view of (4.24)–(4.33), (4.34) entails that

$$\lim_{n \to \infty} \sup_{h \in [h'_n, h''_n]} \frac{|\sqrt{n}\, d_n^{(1)}(\lceil nh \rceil)|}{\sqrt{2h \log(1/h)}} = 1 \text{ a.s.} . \tag{4.35}$$

We conclude by routine analysis that, under $H.(1\text{-}2\text{-}3)$ and the assumptions above, (4.35) holds with the formal replacement of $d_n^{(1)}$ by $d_n$. This completes the proof of Theorem 3.2.

### 4.4. Proof of Theorem 3.1

Recall the definition (2.1) of $\beta_n$. Following the lines of the above-given proof of Theorem 3.2, we select $0 \le s, t \le 1$ and set $i = \lceil nt \rceil$ and $j = \lceil ns \rceil$. We then write

$$\beta_n(t) - \beta_n(s) = n^{1/2}\Big(\big(U_{(j)} - U_{(i)} - \frac{j-i}{n}\big) + \frac{j-i}{n} - (t-s)\Big), \qquad (4.36)$$

and observe that, in (4.36),

$$\Big|\frac{j-i}{n} - (t-s)\Big| \le \frac{1}{n}. \qquad (4.37)$$

In view of (4.36) and (4.37), we obtain the proof of (3.3) by similar arguments as in the proof of Theorem 3.2.

### 4.5. Proof of Theorem 3.3

The Parzen-Rosenblatt kernel estimator of the density function $f$ (see, e.g., Parzen (22) and Rosenblatt (24)) is defined, for some kernel $K$ fulfilling $(K.A\text{-}B)$, and a positive bandwidth $h$, by

$$\widetilde{f}_{n,h}(x) := \frac{1}{nh}\sum_{i=1}^n K\Big(\frac{x - X_i}{h}\Big).$$

Recall the notation $J = (u_1, u_2)$. An application of Theorem 1.1 of Varron (31) (see Facts 5.1 and 5.2 in the Appendix) yields readily the following proposition whose proof is omitted. This proposition provides a uniform-in-bandwidth version of Corollary 4 in Einmahl and Mason (14).

**Proposition 4.2.** *Denote by $\mathcal{J}$ a sub-interval of $J$ with non-empty interior. Let $\{h'_n\}_{n\ge 1}$ and $\{h''_n\}_{n\ge 1}$ be two non-random sequences fulfilling the conditions $(H.1\text{-}2\text{-}3)$, with $0 < h'_n \le h''_n < 1$. Then, under $(F.1\text{-}2\text{-}3)$ and $(K.A\text{-}B)$, we have, almost surely,*

$$\lim_{n\to\infty} \sup_{h\in[h'_n, h''_n]} \sup_{x\in\mathcal{J}} \frac{\sqrt{nh}|\widetilde{f}_{n,h}(x) - \mathbb{E}\widetilde{f}_{n,h}(x)|}{\sqrt{2f(x)\log(1/h)}} = \Big\{\int_{\mathbb{R}} K^2(t)dt\Big\}^{1/2}. \qquad (4.38)$$

To complete the proof of Theorem 3.3, we select two sequences $\{h'_n\}_{n\ge 1}$ and $\{h''_n\}_{n\ge 1}$ fulfilling the conditions $(H.1\text{-}2\text{-}3\text{-}4)$, and set $k'_n = nh'_n$ and $k''_n = nh''_n$. Given this notation, we claim that, almost surely as $n \to \infty$,

$$\sup_{x\in[u_{1,n}, u_{2,n}]} \sup_{k\in[k'_n, k''_n]} \Big|R_k(x) - \frac{k}{nf(x)}\Big| \to 0. \qquad (4.39)$$

To see how (4.39) follows from Theorem 3.2, select $x \in [u_{1,n}, u_{2,n}]$ and let $j = j_k(x)$ be the smallest integer for which $X_{(j)} \geq x - R_k(x)/2$. Then, at least one of the relations

$$X_{(j)} = x - \frac{R_k(x)}{2}, \qquad X_{(j+\lfloor k \rfloor - 1)} = x + \frac{R_k(x)}{2}$$

holds. This naturally implies that

$$R_k(x) \leq \max\left\{ [X_{(j+\lfloor k \rfloor)} - X_{(j)}], [X_{(j+\lfloor k \rfloor - 1)} - X_{(j-1)}] \right\},$$

and

$$R_k(x) \geq [X_{(j+\lfloor k \rfloor - 1)} - X_{(j)}].$$

Moreover, observe that, for all $x' \in [X_{(j+\lfloor k \rfloor - 1)}, X_{(j+\lfloor k \rfloor + 1)}]$, $f(x') = f(x) + o(1)$ almost surely as $n \to \infty$, for all $k \in [k'_n, k''_n]$. Putting all these results together, we can apply Theorem 3.2 and conclude to (4.39).

Now, introduce the modified sequences defined by

$$\tilde{h}'_n := \frac{h'_n}{2} \inf_{x \in J} \frac{1}{f(x)} \quad \text{and} \quad \tilde{h}''_n := 2h''_n \sup_{x \in J} \frac{1}{f(x)}. \tag{4.40}$$

We infer from (4.39) and (4.40) that, almost surely as $n \to \infty$,

$$\tilde{h}'_n \leq \inf_{x \in [u_{1,n}, u_{2,n}]} \inf_{k \in [k'_n, k''_n]} R_k(x) \leq \sup_{x \in [u_{1,n}, u_{2,n}]} \sup_{k \in [k'_n, k''_n]} R_k(x) \leq \tilde{h}''_n. \tag{4.41}$$

Recalling from (F.1-2) that $f$ is bounded and strictly positive on $J$, with $J \supset [u_{1,n}, u_{2,n}]$ for all $n \geq 1$, we see that $\{\tilde{h}'_n\}_{n \geq 1}$ [resp. $\{\tilde{h}''_n\}_{n \geq 1}$] fulfills (H.1-2-3-4), whenever such is the case for $\{h'_n\}_{n \geq 1}$ [resp. $\{h''_n\}_{n \geq 1}$]. Moreover, setting $\tilde{u}_1 = Q(t_1)$ and $\tilde{u}_2 = Q(t_2)$, we have $u_{\ell,n} \to u_\ell$, for $\ell = 1, 2$ almost surely as $n \to \infty$ and $\tilde{u}_1 > u_1$ and $\tilde{u}_2 < u_2$ (by (F.1-2)). Thus, in view of Proposition 4.2, it is straightforward that, under the hypotheses of Theorem 3.3,

$$\lim_{n \to \infty} \sup_{x \in [u_{1,n}, u_{2,n}]} \sup_{k \in [k'_n, k''_n]} \frac{\sqrt{k}\, |\widehat{f}_{n,k}(x) - \mathbb{E}\widehat{f}_{n,k}(x)|}{\sqrt{2f^2(x) \log(n/k)}} \leq \left\{ \int_{\mathbb{R}} K^2(t) dt \right\}^{1/2}. \tag{4.42}$$

Moreover, in view of (4.39), for all $k_n = nh_n \in [k'_n, k''_n]$, it holds that, almost surely as $n \to \infty$,

$$\inf_{x \in J} \frac{1}{f(x)} \frac{h_n}{2} \leq \inf_{x \in [u_{1,n}, u_{2,n}]} R_{k_n}(x) \leq \sup_{x \in [u_{1,n}, u_{2,n}]} R_{k_n}(x) \leq 2h_n \sup_{x \in J} \frac{1}{f(x)}. \tag{4.43}$$

We now compare (4.43) with the condition $(B.1)$ in Deheuvels and Mason (13). The fact that

$$\lim_{n \to \infty} \sup_{x \in [u_{1,n}, u_{2,n}]} \frac{\sqrt{k_n}\, |\widehat{f}_{n,k_n}(x) - \mathbb{E}\widehat{f}_{n,k_n}(x)|}{\sqrt{2f^2(x) \log(n/k_n)}} = \left\{ \int_{\mathbb{R}} K^2(t) dt \right\}^{1/2} \tag{4.44}$$

follows from Proposition 4.2, in view of (4.39) and (4.41), along the same lines as the first part of Theorem 1.2 in Deheuvels and Mason (13) is shown to be a consequence of their Corollary 3.2. We omit the details of this book-keeping argument.

By combining (4.42) with (4.44), the proof of Theorem (3.3) is readily achieved.

**Acknowledgments:** We thank the referees for careful reading of our manuscript and for insightful comments leading, in particular, to some extensions of our original Theorem 1.1.

## 5. Appendix

### 5.1. Useful facts

Here we present two necessary facts, that have been shown to be instrumental in our proofs.

Fix $d \geq 1$. Let $(\mathbf{Z}_i)_{i \geq 1}$ be a sequence of iid random variables defined on $\mathbb{R}^d$, and $\mathcal{G}$ be a class of real Borel functions defined on $\mathbb{R}^d$. For each $h > 0$, $n \geq 1$ and $\mathbf{z} \in \mathbb{R}^d$, we set

$$G_n(K, h, \mathbf{z}) := \sum_{i=1}^n K\left(\frac{\mathbf{z} - \mathbf{Z}_i}{h^{1/d}}\right) - \mathbb{E}\left\{K\left(\frac{\mathbf{z} - \mathbf{Z}_i}{h^{1/d}}\right)\right\},$$

where $K \in \mathcal{G}$.

Further set $I^d := [0,1]^d$, $|u| := \max_{1 \leq i \leq d} |u_i|$ and

$$\mathcal{F} := \left\{K(\lambda(\mathbf{z} - .)), \mathbf{z} \in \mathbb{R}^d, \lambda > 0, K \in \mathcal{G}\right\}.$$

Introduce the following assumptions on $\mathcal{G}$.

$(HK.I)$ $(i)$ $\lim_{|\mathbf{u}| \to \mathbf{0}} \sup_{K \in \mathcal{G}} \int_{\mathbb{R}^d} \left(K(\mathbf{x}) - K(\mathbf{x} + \mathbf{u})\right)^2 d\mathbf{x} = 0$.

$(ii)$ $\lim_{\lambda \to 1} \sup_{K \in \mathcal{G}} \int_{\mathbb{R}^d} \left(K(\lambda \mathbf{x}) - K(\mathbf{x})\right)^2 d\mathbf{x} = 0$.

$(HK.II)$ $(i)$ $\forall K \in \mathcal{G}$, $\sup_{\mathbf{x} \in \mathbb{R}^d} |K(\mathbf{x})| \leq 1$.

$(ii)$ $\forall K \in \mathcal{G}$, $\forall x \notin I^d$, $K(\mathbf{x}) = 0$.

$(HK.III)$ $(i)$ $\exists\, C > 0, v > 0, \forall\, 0 < \varepsilon < 1$, $\mathcal{N}(\varepsilon, \mathcal{F}) \leq C\varepsilon^{-v}$.

$(ii)$ $\mathcal{F}$ is pointwise separable.

Here $\mathcal{N}(\varepsilon, \mathcal{F}) := \sup\{\mathcal{N}(\varepsilon, \mathcal{F}, \mathbb{L}_2(\mathbb{P})), \mathbb{P} \text{ probability measure}\}$ denotes the uniform covering number of $\mathcal{F}$ for $\varepsilon$ and the class of norms $\{\mathbb{L}_2(\mathbb{P})\}$, with $\mathbb{P}$ varying in the set of all probability measures on $\mathbb{R}^d$ (for more details we refer to van der Vaart and Wellner (28), pp 83-84).

Let $\mathcal{B}(\mathcal{G})$ denote the set of all real bounded functions on $\mathcal{G}$, continuous with respect to pointwise convergence. Let $\mathbb{L}_2^\star(\mathcal{G})$ be the Hilbert subspace of $\mathbb{L}_2(\mathbb{R}^d, m)$ spanned by $\mathcal{G}$ (here, $m$ denotes the Lebesgue measure on $\mathbb{R}^d$). For $f \in \mathbb{L}_2^\star(\mathbb{R}^d) = \mathbb{L}_2(\mathbb{R}^d)$, denote by $\Psi_f(g) := (f,g) = \sqrt{\int fg dm}$, $g \in \mathcal{G}$, the inner product in $\mathbb{L}_2(\mathbb{R}^d, m)$ and $J(\Psi_f) := \int_{\mathbb{R}^d} f^2 dm$ the squared norm of $f$ in $\mathbb{L}_2(\mathbb{R}^d, m)$. For each $\Psi \in \mathcal{B}(\mathcal{G})$ that can not be expressed as a $\Psi_f$, $f \in \mathbb{L}_2(\mathbb{R}^d)$, we set $J(\Psi) = \infty$. Finally, set $\mathbb{K} := \{\Psi : J(\Psi) \leq 1\}$.

The following result is due to Varron (31) (see also Varron (32)).

**Fact 5.1.** *Let $\{h'_n\}_{n \geq 1}$ and $\{h''_n\}_{n \geq 1}$ be two sequences of positive constants, fulfilling the assumptions (H.1-2-3), and such that $0 < h'_n \leq h''_n < 1$. Suppose that $\mathbf{Z}_1$ has a density function $f$, such that the following conditions hold. For some compact set $H \in \mathbb{R}^d$, there exists a $\varsigma > 0$ such that $f$ is continuous and strictly positive on*

$$H^\varsigma := \{\boldsymbol{x} \in \mathbb{R}^d : \inf_{\mathbf{z} \in H} \|\boldsymbol{x} - \mathbf{z}\|_{\mathbb{R}^d} \leq \varsigma\}.$$

*Here, $\|.\|_{\mathbb{R}^d}$ denotes the usual Euclidean norm on $\mathbb{R}^d$. Then, if $\mathcal{G}$ is a class of real Borel functions satisfying $(HK.I\text{-}II\text{-}III)$, we have almost surely,*

$$\lim_{n\to\infty} \sup_{\mathbf{z} \in H} \sup_{h \in [h'_n, h''_n]} \inf_{\Psi \in \mathbb{K}} \left\| \frac{G_n(.,h,\mathbf{z})}{\sqrt{2f(\mathbf{z})nh\log(1/h)}} - \Psi \right\|_{\mathcal{G}} = 0,$$

$$\forall \Psi \in \mathbb{K}, \lim_{n\to\infty} \sup_{h \in [h'_n, h''_n]} \inf_{\mathbf{z} \in H} \left\| \frac{G_n(.,h,\mathbf{z})}{\sqrt{2f(\mathbf{z})nh\log(1/h)}} - \Psi \right\|_{\mathcal{G}} = 0,$$

*where $\|\psi\|_{\mathcal{G}} = \sup_{K \in \mathcal{G}} |\psi(K)|$.*

A direct consequence of this result is as follows (see, e.g., Varron (32)).

**Fact 5.2.** *Let $K$ be a measurable kernel of bounded variation, with compact support. Under the relevant assumptions of Fact 5.1, we have, almost surely,*

$$\lim_{n\to\infty} \sup_{h \in [h'_n, h''_n]} \sup_{\boldsymbol{x} \in H} \frac{\sqrt{nh}|f_n(K, \boldsymbol{x}, h) - \mathbb{E}f_n(K, \boldsymbol{x}, h)|}{\sqrt{2f(\boldsymbol{x})\log(1/h)}} = \sqrt{\int_{\mathbb{R}^d} K^2 dm},$$

*where $f_n(K, \boldsymbol{x}, h) = \dfrac{1}{nh} \sum_{i=1}^{n} K\left(\dfrac{\mathbf{z} - \mathbf{Z}_i}{h^{1/d}}\right) - \mathbb{E}\left\{K\left(\dfrac{\mathbf{z} - \mathbf{Z}_i}{h^{1/d}}\right)\right\}.$*

## References

[1] R.R. Bahadur. A note on quantiles in large samples. *Ann. Math. Stat.*, 37:577–580, 1966. MR0189095
[2] K.L. Chung. On the maximum partial sums of sequences of independent random variables. *Trans. Am. Math. Soc.*, 64:205-233, 1948. MR0026274
[3] M. Csörgő and P. Révész. Strong approximations of the quantile process. *Ann. Stat.*, 6:882–894, 1978. MR0501290

[4] M. Csörgő and P. Révész. *Strong Approximation in Probability and Statistics.* Acedemic Press, New York, 1981. MR0666546
[5] M. Csörgő and P. Révész. Two approaches to constructing simultaneous confidence bounds for quantiles. *Probab. Math. Statist.*, 4(2):221–236, 1984. MR0792787
[6] M. Csörgő. *Quantile processes with statistical applications.* CBMS-NSF Regional Conference Series in Applied Mathematics, SIAM, Philadelphia, 1983. MR0745130
[7] P. Deheuvels. Conditions nécessaires et suffisantes de convergence ponctuelle presque sûre des estimateurs de la densité. *C. R. Acad. Sci., Paris, Sér. A*, 278:1217–1220, 1974. MR0345296
[8] P. Deheuvels. Strong limiting bounds for maximal uniform spacings. *Ann. Probab.*, 10:1058–1065, 1982. MR0672307
[9] P. Deheuvels. Functional laws of the iterated logarithm for large increments of empirical and quantile processes. *Stochastic Processes Appl.*, 43(1):133–163, 1992. MR1190911
[10] P. Deheuvels and J.H.J. Einmahl. Functional limit laws for the increments of Kaplan-Meier product-limit processes and applications. *Ann. Prob*, 28(7):1301–1335, 2000. MR1797314
[11] P. Deheuvels and D.M. Mason. Nonstandard functional laws of the iterated logarithm for tail empirical and quantile processes. *Ann. Prob.*, 18:1693–1722, 1990. MR1071819
[12] P. Deheuvels and D.M. Mason. Functional laws of the iterated logarithm for the increments of empirical and quantile processes. *Ann. Prob.*, 20:1248–1287, 1992. MR1175262
[13] P. Deheuvels and D.M. Mason. General confidence bounds for nonparametric functional estimators. *Stat. Inf. for Stoch. Proc.*, 7:225–277, 2004. MR2111291
[14] U. Einmahl and D.M. Mason. An empirical process approach to the uniform consistency of kernel type estimators. *Journ. Theoretic. Probab.*, 13:1–13, 2000. MR1744994
[15] U. Einmahl and D.M. Mason. Uniform in bandwidth consistency of kernel-type function estimators. *Ann. Stat.*, 33(3):1380–1403, 2005. MR2195639
[16] J. Kiefer. On Bahadur's representation of sample quantiles. *Ann. Math. Stat.*, 38:1323–1342, 1967. MR0217844
[17] J. Kiefer. Iterated logarithm analogues for sample quantiles when $p_n \downarrow 0$. In *Proc. 6th Berkeley Sympos. math. Statist. Probab.*, pages 227–244. Univ. California Press, Berkeley, 1972. MR0402882
[18] D.M. Mason. A strong limit theorem for the oscillation modulus of the uniform empirical quantile process. *Stochastic Process. Appl.*, 17:127–136, 1984. MR0738772
[19] D.M. Mason. A uniform functional law of the iterated logarithm for the local empirical process. *Ann.Prob.*, 32(2):1391–1418, 2004. MR2060302
[20] D.M. Mason, G.R. Shorack, and J.A. Wellner. Strong limit theorems for oscillation moduli of the uniform empirical process. *Z. Wahrscheinlichkeitstheor. Verw. Geb.*, 65:83–97, 1983. MR0717935

[21] E.A. Nadaraya. On estimating regression. *Theor.Prob.Appl.*, 9:141–142, 1964.
[22] E. Parzen. On estimation of probability density function and mode. *Ann.Math.Stat.*, 33:1065–1076, 1962. MR0143282
[23] E. Parzen. Nonparametric statistical data modelling. *J. Amer. Statist. Assoc.*, 74:105–131, 1979. MR0529528
[24] M. Rosenblatt. Remarks on some nonparametric estimates of a density function. *Ann.Math.Stat.*, 27:832–837, 1956. MR0079873
[25] V. Strassen. An invariance principle for the law of the iterated logarithm. *Z.Wahrsch.Gebiete*, 3:221–226, 1964. MR0175194
[26] W. Stute. The oscillation behavior of empirical processes. *Ann. Prob.*, 10:86–107, 1982. MR0637378
[27] W. Stute. The law of the iterated logarithm for kernel density estimators. *Ann. Prob.*, 10:414–422, 1982. MR0647513
[28] A.W. van der Vaart and J.A. Wellner. *Weak convergence and empirical processes.* Springer, New York, 1996. MR1385671
[29] V.N. Vapnik. *Estimation of dependences.* Springer-Verlag, N.Y., 1982.
[30] V.N. Vapnik and A.Y. Červonenkis. On the uniform convergence of relative frequencies of events to their probabilities. *Theor. Probab. Appl.*, 16:264–280, 1971.
[31] D. Varron. Uniformity in $h$ in the functional limit law for the increments of the empirical process indexed by functions. (Uniformité en $h$ dans la loi fonctionnelle limite uniforme des accroissements du processus empirique indexé par des fonctions). *C. R., Math., Acad. Sci. Paris*, 340(6):453–456, 2005. MR2135329
[32] D. Varron. A bandwidth-uniform functional limit law for the increments of the empirical process. *Preprint*, 2006.
[33] G.S. Watson. Smooth regression analysis. *Sankhya, Ser. A*, 26:359–372, 1964. MR0185765